\documentclass[10pt,letterpaper]{amsart}

\usepackage[utf8]{inputenc}
\usepackage{amsmath,amssymb,mathtools}

\usepackage[unicode]{hyperref}
\usepackage{enumitem}

\usepackage{color}

\allowdisplaybreaks[2]

\newtheorem{theorem}{Theorem}[section]
\newtheorem{lemma}[theorem]{Lemma}
\newtheorem{proposition}[theorem]{Proposition}

\theoremstyle{definition}
\newtheorem{definition}[theorem]{Definition}
\theoremstyle{remark}

\newcommand{\N}{\mathbb N}
\newcommand{\PP}{\mathcal P}
\newcommand{\restr}{\mathbin{\upharpoonright}}

\newcommand{\sm}{\mathrm{sm}}
\newcommand{\Lsm}{L_{\mathrm{sm}}}
\newcommand{\Lce}{L^0_{\mathrm{sm}}}
\newcommand{\smle}{\leq_{\mathrm{sm}}}
\newcommand{\smeq}{\equiv_{\mathrm{sm}}}
\newcommand{\dsm}{\operatorname{deg}_{\mathrm{sm}}}

\newcommand{\Th}{\operatorname{Th}}
\newcommand{\TAtwo}{\mathrm{TA}_2}

\newcommand{\Prm}{\operatorname{JI}}
\newcommand{\CP}{\operatorname{Cup}}

\newcommand{\Code}{\operatorname{Code}}

\newcommand{\EqSet}{\operatorname{EqSet}}

\title[Separation degrees]{On the theory of separation degrees}
\date{\today}

\author{Nikolay Bazhenov}
\address{Nanyang Technological University, 21 Nanyang Link, 637371, Singapore}
\email{nikolay.bazhenov@ntu.edu.sg}


\begin{document}

\begin{abstract}
 	Ershov (1970) introduced the \emph{separation reducibility} $\leq_{\mathrm{sm}}$ on ordered pairs of disjoint sets. A pair $(A_0,A_1)$ is $\mathrm{sm}$-reducible to a pair $(B_0,B_1)$ iff there exists a total computable function $f$ such that $x\in A_i$ implies $f(x) \in B_i$. Dvornikov and Rybina (1990) investigated the structural properties of the poset $L_{\mathrm{sm}}$ of all separation degrees. In particular, they showed that $L_{\mathrm{sm}}$ is a dense distributive lattice. The question of the complexity of the first-order theory of $L_{\mathrm{sm}}$ was left open in previous studies. In this paper we prove that the theory $\operatorname{Th}(L_{\mathrm{sm}})$ is recursively isomorphic to the second-order arithmetic.
\end{abstract}

\maketitle

\section{Introduction}

Let $A$ and $B$ be subsets of the set of natural numbers $\omega$. We say that the (ordered) pair $(A,B)$ is \emph{disjoint} if $A\cap B = \emptyset$. In his studies of the theory of numberings, Ershov~\cite{Ershov70} (see also Chapter~3.2 in~\cite{Ershov-Book}) introduced the notion of the \emph{separation reducibility} $\smle$: for disjoint pairs $(A_0,B_0)$ and $(A_1,B_1)$, one says that $(A_0,B_0) \smle (A_1,B_1)$ iff there exists a total computable function $f(x)$ satisfying the following two conditions:
\begin{itemize}
	\item if $x\in A_0$, then $f(x) \in A_1$;
	
	\item if $x\in B_0$, then $f(x) \in B_1$.
\end{itemize}
As usual, we say that $(A_0,B_0)\smeq (A_1,B_1)$ iff $(A_0,B_0) \smle (A_1,B_1)$ and $(A_1,B_1) \smle (A_0,B_0)$.

The resulting poset $\Lsm$ containing separation degrees of all disjoint pairs exhibits quite a curious behavior. In contrast to many familiar degree structures in computability theory (e.g., many-one degrees, truth-table degrees, or Turing degrees), the poset $\Lsm$ is a \emph{lattice}, where:
\begin{itemize}
	\item as usual, the join of $(A,B)$ and $(C,D)$ is defined as $(A\oplus C, B\oplus D)$;
	
	\item the meet is given as $(A\times C, B\times D)$, where $A\times C = \{ \langle a,c\rangle : a\in A,\, c\in C\}$ (see \cite{DvR90}).
\end{itemize}
In addition, the lattice $\Lsm$ is distributive and dense (Theorem~1 and Corollary~1 of~\cite{DvR90}). 

The papers~\cite{DvR90,Ryb90} continued investigations of the separation reducibility. In particular, Rybina~\cite{Ryb90} considered the lattice $\Lce$ containing separation degrees of computably enumerable (c.e.) disjoint pairs, i.e., pairs $(A,B)$ such that both sets $A$ and $B$ are c.e. She proved that the first-order theory $\Th(\Lce)$ is hereditarily undecidable.

Dvornikov and Rybina~\cite{DvR90} left open the questions about the exact complexity of the theories $\Th(\Lce)$ and $\Th(\Lsm)$. In this paper we prove the following: 

\begin{theorem}\label{thm:main}
	For the lattice $\Lsm$ of all separation degrees, its first-order theory $\Th(\Lsm)$ is recursively isomorphic to the second-order arithmetic.
\end{theorem}

The paper is arranged as follows. Section~\ref{sec:prelim} contains the necessary preliminaries, and it outlines the proof of Theorem~\ref{thm:main}. The rest of the paper gives a detailed proof of Theorem~\ref{thm:main}. We leave open whether the theory $\Th(\Lce)$ of c.e.\ separation degrees is recursively isomorphic to the first-order arithmetic.

\medskip

\section{Preliminaries}\label{sec:prelim}

We refer to~\cite{Odi-92,Soare-87} for the background on computability theory. Following~\cite{Ershov-Book,DvR90}, we assume that $\mathbf{0} = \dsm (\{0\}, \{1\})$ is the least degree in the lattice $\Lsm$. Notice that every recursively separable pair of non-empty sets $(A,B)$ belongs to $\mathbf{0}$.

Since $\Lsm$ is a lattice with a least element, we will assume that the language of the structure $\Lsm$ is equal to $\{ \leq, \vee,\wedge,\mathbf{0}\}$.

In what follows, we denote disjoint pairs by underlined capital letters: $\underline{A},\underline{B},\underline{C},\dots$. We always assume that $\varepsilon \in \{ 0,1\}$. For a disjoint pair $\underline{A}$, by $\underline{A}^{[\varepsilon]}$ we denote the $\varepsilon$-th component of $\underline{A}$---that is, $\underline{A} = (\underline{A}^{[0]}, \underline{A}^{[1]})$.

For a countable sequence of disjoint pairs $(\underline{A_i})_{i\in\omega}$, its \emph{uniform join} is the disjoint pair $\underline{B} = \bigoplus_{i\in\omega} \underline{A_i}$ such that
\begin{equation*} 
	\underline{B}^{[\varepsilon]} = \{\langle i,x\rangle : x \in \underline{A_i}^{[\varepsilon]},\, i\in\omega\}.
\end{equation*}

\medskip

For a computable set $R\subseteq \omega$, we define the \emph{$R$-restriction} of a disjoint pair $\underline{A}$ as follows:
\begin{equation*} 
 (\underline{A\restr R})^{[\varepsilon]} = \{\varepsilon\} \cup \{ 2x+2 : x\in \underline{A}^{[\varepsilon]}\cap R\}.
\end{equation*}

If $\mathbf{a} = \dsm(\underline{A})$, then we write
\[
	\mathbf{a}_R = \dsm(\underline{A\restr R}). 
\]
Note that the notation $\mathbf{a}_R$ depends on the chosen representative $\underline{A}$ in $\mathbf{a}$, but in later proofs $\underline{A}$ will remain fixed whenever $R$-restrictions are considered. Observe that
\begin{equation}\label{eq:split}
 \mathbf{a}_R\leq \mathbf{a},\quad
  \mathbf{a} = \mathbf{a}_R\vee \mathbf{a}_{\overline R}.
\end{equation}
It is also easy to show the following:

\begin{lemma}\label{lem:tagcut}
	Suppose that $\underline{A}\smle \underline{P}\oplus \underline{Q}$ via a function $f$. Let $R=\{x : f(x) \text{ is even}\}$. Then we have:
	\[
		\underline{A\restr R} \smle \underline{P},\quad
		\underline{A\restr\overline R} \smle \underline{Q}.
	\]
\end{lemma}

\medskip

We will use the following result of Dvornikov and Rybina:

\begin{theorem}[Theorem~4 of~\cite{DvR90}]\label{thm:separation}
	Suppose that $(\mathbf{a}_n)_{n\in\omega}$ and $(\mathbf{d}_n)_{n\in\omega}$ are sequences of $\sm$-degrees such that
	\[
		 \mathbf{a}_n < \mathbf{a}_{n+1} \text{ and } \mathbf{d}_n\nleq \mathbf{a}_m,\quad \text{for all }n,m\in\omega.
	\]
	Then there exists a degree $\mathbf{c}$ such that
	\[
 		\mathbf{a}_n\leq \mathbf{c} \text{ and } \mathbf{d}_n\nleq \mathbf{c},\quad \text{for all } n\in\omega.
	\]
\end{theorem}

\medskip

\subsection{Outline of the proof of the main result} 

Our main goal is to `encode' the second-order arithmetic
\[
	\TAtwo = \Th\bigl(\N,\PP(\N);+,\times,\in\bigr)
\]
into the theory $\Th(\Lsm)$. 

First, in Section~\ref{sec:sets} we show how to give a first-order interpretation of the powerset structure $(\N,\PP(\N);\in)$ inside $\Lsm$, with one parameter $\mathbf{a}$.
In Section~\ref{sec:free} we explicitly construct such parameters $\mathbf{a}$, requiring some additional properties of $\mathbf{a}$.
In Section~\ref{sec:relations} we show how to encode a given countable undirected graph into $\Lsm$ using parameters.
Finally, in Section~\ref{sec:finish} we finish the proof of Theorem~\ref{thm:main}: 
\begin{enumerate}
	\item[(i)] we fix a computable undirected graph $G$ such that the structure $(\N,+,\times)$ is first-order definable (without parameters);
	
	\item[(ii)] using the constructions of Sections~\ref{sec:free} and~\ref{sec:relations} we simultaneously encode $(\N;+,\times)$ and $(\N,\PP(\N);\in)$ inside $\Lsm$ using parameters $\mathbf{a},\mathbf{b},\mathbf{g}$;
	
	\item[(iii)] we show how to eliminate the parameters $\mathbf{a},\mathbf{b},\mathbf{g}$ from Item~(ii), so our encoding gives the desired effective encoding of $\TAtwo$ into $\Th(\Lsm)$.
\end{enumerate}

\medskip

\section{Interpreting the powerset inside $\Lsm$}\label{sec:sets}

In this section we provide a first-order interpretation of the (two-sorted) powerset structure $(\N,\PP(\N);\in)$ inside $\Lsm$, with one parameter $\mathbf{a}$. We leave the proof of the existence of an appropriate parameter $\mathbf{a}\in \Lsm$ until the next section.

First, we give the necessary preliminary definitions. 

\begin{definition}\label{def:Prm}
	We introduce the formula
	\begin{equation*} 
		\Prm(x) := (x\ne\mathbf{0})\ \&\  \forall u\, \forall v\, \bigl(x\leq u\vee v\ \longrightarrow\ x\leq u \text{ or } x\leq v\bigr).
	\end{equation*}
\end{definition}

Note that in the lattice $\Lsm$, the condition $\Prm(x)$ is equivalent to $x$ being non-zero and join-irreducible. Indeed, in any lattice $L$, $\Prm(x)$ implies join-irreducibility. Now assume that $x\neq \mathbf{0}$ is join-irreducible. If $x \leq u\vee v$, then by the distributivity of $\Lsm$ we obtain 
\[
	x= x\wedge (u\vee v) = (x\wedge u)\vee(x\wedge v),
\]
and the join-irreducibility of $x$ implies that $x = x\wedge u$ or $x = x\wedge v$.

\smallskip

\begin{definition}\label{def:Ja}
	For a non-zero degree $\textbf{a}\in\Lsm$ we define the formula
	\begin{equation*} 
		\CP_{\mathbf{a}}(x) := \Prm(x)\ \&\ x<\mathbf{a}\ \&\ \exists y\,(y<\mathbf{a}\ \&\ \mathbf{a}=x\vee y).
	\end{equation*}
	That is, $\CP_{\mathbf{a}}(x)$ holds if and only if $x$ is a non-zero join-irreducible element which is cuppable to $\mathbf{a}$.
\end{definition}

We show the following:

\begin{lemma}\label{lem:Ja}
	For every non-zero $\sm$-degree $\mathbf{a}$, the set $\CP_{\mathbf{a}} [\Lsm] = \{ x : \CP_{\mathbf{a}}(x)\}$ is an at most countable antichain.
\end{lemma}
\begin{proof}
Towards a contradiction, assume that degrees $\mathbf{x} < \mathbf{z}$ belong to $\CP_{\mathbf{a}} [\Lsm]$. 
We choose $\mathbf{y}<\mathbf{a}$ such that $\mathbf{a} = \mathbf{x} \vee \mathbf{y}$. Since $\mathbf{z}\leq \mathbf{a}$ satisfies $\Prm(\mathbf{z})$ and $\mathbf{z}\nleq \mathbf{x}$, we have $\mathbf{z}\leq \mathbf{y}$. Then $\mathbf{x} < \mathbf{z}\leq \mathbf{y}$ and $\mathbf{a}= \mathbf{x} \vee \mathbf{y} = \mathbf{y} <\mathbf{a}$, a contradiction. We deduce that the set $\CP_{\mathbf{a}}[\Lsm]$ is an antichain.

Now fix a pair $\underline{A} \in \mathbf{a}$. Given an arbitrary $\mathbf{x}\in \CP_{\mathbf{a}}[\Lsm]$, we choose $\mathbf{y}<\mathbf{a}$ satisfying $\mathbf{a}=\mathbf{x}\vee \mathbf{y}$. We also choose pairs $\underline{X}$ and $\underline{Y}$ from $\mathbf{x}$ and $\mathbf{y}$, respectively.
For a fixed computable reduction $f\colon \underline{A}\smle \underline{X}\oplus \underline{Y}$ define the set $R = \{ z : f(z) \text{ is even}\}$. By Eq.~(\ref{eq:split}) and Lemma~\ref{lem:tagcut}, we have:
\[
 	\mathbf{a} = \mathbf{a}_R\vee \mathbf{a}_{\overline R},\quad \mathbf{a}_R\leq \mathbf{x},\quad
	 \mathbf{a}_{\overline R}\leq \mathbf{y}.
\]
Since $\mathbf{x}\leq \mathbf{a}$, $\Prm(\mathbf{x})$ implies that $\mathbf{x}\leq \mathbf{a}_R$ or $\mathbf{x}\leq \mathbf{a}_{\overline R}$. If $\mathbf{x}\leq \mathbf{a}_{\overline R}$, then $\mathbf{x}\leq \mathbf{y}$ and $\mathbf{a}=\mathbf{y}$, a contradiction. 
Hence, we have $\mathbf{x} = \mathbf{a}_R$. Since there are only countably many computable sets $R$, there are at most countably many elements $\mathbf{x}\in \CP_{\mathbf{a}}[\Lsm]$.
\end{proof}

\smallskip

\begin{definition}\label{def:join-ind}
	We say that a family of non-zero $\sm$-degrees $T$ is \emph{join-independent} iff
	\begin{equation}\label{eq:independent}
		 \mathbf{t}\nleq\bigvee F \text{ for every }\mathbf{t}\in T\text{ and every non-empty finite }F\subseteq T\setminus\{\mathbf{t}\}.
	\end{equation}
\end{definition}

We show the following:

\begin{lemma}\label{lem:trace}
	Let $T$ be an at most countable join-independent family with an upper bound $\mathbf{a}$. Then for every $S\subseteq T$ there exists $\mathbf{c}\leq \mathbf{a}$ such that
	\begin{equation}\label{eq:trace}
		 \forall \mathbf{t}\in T\ (\mathbf{t}\leq \mathbf{c}\ \leftrightarrow\ \mathbf{t}\in S).
	\end{equation}
\end{lemma}
\begin{proof}
If $S=\emptyset$, then one can choose $\mathbf{c}=\mathbf{0}$. If $S=T$, then take $\mathbf{c}=\mathbf{a}$. If $S$ is finite and nonempty, then we take its join $\mathbf{c} = \bigvee S$, and the join-independence of $T$ (see Eq.~(\ref{eq:independent})) immediately implies Eq.~(\ref{eq:trace}). 

Suppose that $S\neq T$ and $S = \{\mathbf{s}_i : i\in\omega\}$, where $\mathbf{s}_i \neq \mathbf{s}_j$ for $i\neq j$. We define $\mathbf{a}_n=\bigvee_{i\leq n}\mathbf{s}_i$. Since $T$ is join-independent, we have $\mathbf{s}_{n+1} \nleq \mathbf{a}_n$ and hence, $\mathbf{a}_n < \mathbf{a}_{n+1}$ for all $n$.

Since $T\setminus S$ is nonempty, we assume that $T\setminus S = \{\mathbf{d}_i : i\in\omega\}$, where, in general, $\mathbf{d}_i$ and $\mathbf{d}_j$ could be equal for $i\neq j$. The join-independence of $T$ implies that $\mathbf{d}_n\nleq \mathbf{a}_m$ for all $n,m\in \omega$. We apply Theorem~\ref{thm:separation} to the sequences $(\mathbf{a}_i)_{i\in\omega}$ and $(\mathbf{d}_n)_{n\in\omega}$. This gives us an upper bound $\mathbf{c}'$ of the degrees $\mathbf{a}_n$, $n\in\omega$, such that every $\mathbf{d}_n$ satisfies $\mathbf{d}_n\nleq \mathbf{c}'$. Now our desired degree $\mathbf{c}$ is defined as $\mathbf{c} =\mathbf{c}'\wedge \mathbf{a}$. Since all members of $T$ are below $\mathbf{a}$, we still have $\mathbf{s}_n \leq \mathbf{a}_n \leq \mathbf{c}$ and $\mathbf{d}_n\nleq \mathbf{c}$. 
Lemma~\ref{lem:trace} is proved.
\end{proof}

\medskip

We obtain the main result of this section.

\begin{proposition}\label{prop:fullsets}
	Let $\mathbf{a}$ be a non-zero $sm$-degree.
	\begin{itemize}
		\item[(a)] For every $D\subseteq \CP_{\mathbf{a}}[\Lsm]$ and every $S\subseteq D$, there exists $\mathbf{c}\leq \mathbf{a}$ satisfying
	\begin{equation*} 
		 S=\{\mathbf{x}\in D : \mathbf{x} \leq \mathbf{c}\}.
	\end{equation*}
	
		\item[(b)] If the set $\CP_{\mathbf{a}}[\Lsm]$ is infinite, then the structure $(\N,\PP(\N);\in)$ has a first-order interpretation with one parameter $\mathbf{a}$ (inside $\Lsm$).
	\end{itemize}
\end{proposition}
\begin{proof}
	\textbf{(a)}\ The set $\CP_{\mathbf{a}}[\Lsm]$ is join-independent. Indeed, assume that $\mathbf{x} \leq \bigvee F$ for a degree $\mathbf{x} \in \CP_{\mathbf{a}}[\Lsm]$ and finite non-empty $F \subseteq \CP_{\mathbf{a}}[\Lsm]$. Then, $\Prm(\mathbf{x})$ implies that $\mathbf{x} \leq \mathbf{y}$ for some $\mathbf{y}\in F$. By Lemma~\ref{lem:Ja}, $\CP_{\mathbf{a}}[\Lsm]$ is an antichain, hence, $\mathbf{x} = \mathbf{y}$.
	
	Now we can apply Lemma~\ref{lem:trace} to the set $D\subseteq \CP_{\mathbf{a}}[\Lsm]$ and its upper bound $\mathbf{a}$.
	
	\smallskip
	
	\textbf{(b)}\ Here by a definable relation we mean a relation first-order definable with one parameter $\mathbf{a}$ inside the lattice $\Lsm$. Suppose that the set $\CP_{\mathbf{a}}[\Lsm]$ is infinite. Note that by Lemma~\ref{lem:Ja}, $\CP_{\mathbf{a}}[\Lsm]$ is countable. Then we can choose a countably infinite definable subset $D\subseteq \CP_{\mathbf{a}}[\Lsm]$ (say, just put $D = \CP_{\mathbf{a}}[\Lsm]$). 
	
	We use the following first-order interpretation of the two-sorted structure $(\N,\PP(\N);\in)$:
	\begin{itemize}
		\item An element $x\in \N$ is interpreted as a degree $\mathbf{x}\in D$.
		
		\item An element $C\in\PP(\N)$ is interpreted as its code $\mathbf{c} \leq \mathbf{a}$.
		
		\item For $x\in\N$ and $C\in \PP(\N)$, the condition $x\in C$ is interpreted as $\mathbf{x} \leq \mathbf{c}$. 
	\end{itemize}
	Note that here Item~(a) guarantees that every set $C\subseteq \N$ obtains at least one code.
	In general, one set $C$ can get different codes. So, we have to take a quotient of the constructed model $\mathcal{M} = (D, \{\mathbf{c} : \mathbf{c}\leq \mathbf{a}\};\leq)$ via the definable relation
	\begin{equation*}
		 \EqSet_D(u,v) := \forall x\,\bigl(D(x)\rightarrow(x\leq u\leftrightarrow x\leq v)\bigr).
	\end{equation*}
	(Here the quotient is taken only on the second sort of the model, i.e., on $\{\mathbf{c} : \mathbf{c}\leq \mathbf{a}\}$). It is clear that the quotient structure $\mathcal{M}/\EqSet_D$ is isomorphic to $(\N,\PP(\N);\in)$. Proposition~\ref{prop:fullsets} is proved.
\end{proof}

\medskip

\section{Antichains of join-irreducible degrees}\label{sec:free}

In this section we construct some specific $\sm$-degrees $\mathbf{a}$ such that the set $\CP_{\mathbf{a}}[\Lsm]$ is a countably infinite antichain (recall Lemma~\ref{lem:Ja}). Note that we also require some additional properties of $\mathbf{a}$ (which will be used later). Applying Proposition~\ref{prop:fullsets} to such $\mathbf{a}$ allows us to interpret the powerset structure $(\N, \PP(\N);\in)$ inside the down-set $(\{ \mathbf{c} : \mathbf{c} \leq \mathbf{a}\}, \smle)$.

Recall that an infinite set $C\subseteq\omega$ is \emph{$r$-cohesive} if, for every computable set $R\subseteq \omega$, either $C\cap R$ or $C\cap \overline{R}$ is finite. Note that $r$-cohesive sets exist (see, e.g., Proposition~III.4.16 in~\cite{Odi-92}).

First, we establish the following ancillary fact about $r$-cohesive sets.

\begin{lemma}\label{lem:cohesiveprime}
	Suppose that a set $C$ is $r$-cohesive, and $C$ is a disjoint union of two infinite sets $A$ and $B$. Then the degree $\mathbf{x} = \dsm(A,B)$ satisfies the formula $\Prm(x)$ from Definition~\ref{def:Prm}--- that is, $\mathbf{x}$ is nonzero and join-irreducible.
\end{lemma}
\begin{proof}
Assume that $\mathbf{x} = \mathbf{0}$. Then there is a computable set $R$ such that $A\subseteq R$ and $B\subseteq \overline{R}$. Hence, both sets $C\cap R = A$ and $C\cap \overline{R} = B$ are infinite, which contradicts the $r$-cohesiveness of $C$. Therefore $\mathbf{x}$ is non-zero.

Now suppose that $f\colon (A,B) \smle \underline{U} \oplus \underline{V}$. Consider a computable set $S = \{x : f(x) \text{ is even}\}$. Since $C$ is $r$-cohesive, one of the sets $C\cap S$ or $C\cap \overline{S}$ is finite. Without loss of generality, assume that $F = C\cap \overline{S}$ is finite. Without loss of generality, also assume that $0\in \underline{U}^{[0]}$ and $1\in \underline{U}^{[1]}$. Then a computable function
\[
	g(x) = \begin{cases}
		m, & \text{if } f(x) = 2m,\\
		n, & \text{if } f(x) = 2n+1 \text{ and } x\not\in F,\\
		0, & \text{if } x\in F\cap A,\\
		1, & \text{if } x\in F\cap B.
	\end{cases}
\]
provides an $\sm$-reduction from $(A,B)$ to the pair $\underline{U}$. We conclude that the degree $\mathbf{x}$ satisfies $\Prm(\mathbf{x})$. Lemma~\ref{lem:cohesiveprime} is proved.
\end{proof}

\medskip

Suppose that $(\mathbf{p}_i)_{i\in\omega}$ is a sequence of $\sm$-degrees.
For a nonempty finite set $F\subseteq\omega$, we use notation
\[
	\mathbf{p}_F = \bigwedge_{i\in F} \mathbf{p}_i. 
\]
We establish the following main result of this section:

\begin{proposition}\label{prop:free}
	There exist non-zero join-irreducible degrees $\mathbf{p}_i = \dsm(\underline{P_i})$, $i\in\omega$, such that:
	\begin{itemize}
		\item[(a)] for every nonempty finite set $F\subseteq \omega$, the degree $\mathbf{p}_F$ is non-zero;

		\item[(b)] for every nonempty finite $F\subseteq \omega$ and every nonempty finite family $\mathcal{G}$ of nonempty finite subsets of $\omega$, we have
		\begin{equation}\label{eq:free}
 			\mathbf{p}_F \leq \bigvee_{G\in\mathcal{G}}\mathbf{p}_G \ \Longleftrightarrow\ \exists G\in\mathcal{G}\;(G\subseteq F);
		\end{equation}

		\item[(c)] for every $i\in\omega$,
		\begin{equation}\label{eq:omit}
			\underline{P_i} \not\smle \bigoplus_{j\ne i} \underline{P_j},
		\end{equation}
		where $\bigoplus_{j\ne i}\underline{P_j}$ is the uniform join of $\{ \underline{P_j} : j\neq i\}$.
	\end{itemize}
	Consequently, the degrees $\mathbf{p}_i$, together with $\mathbf{0}$, generate a free distributive lattice with a least element. In addition, for the degree $\mathbf{a} = \dsm\left(\bigoplus_{i\in\omega} \underline{P_i}\right)$, the set $\CP_{\mathbf{a}}[\Lsm]$ contains all $\mathbf{p}_i$, $i\in\omega$.
\end{proposition}
\begin{proof}
We fix an $r$-cohesive set $C$. Without loss of generality, we may assume that the set $C$ is co-c.e.\ (consequently, $C\leq_T \mathbf{0}'$). We will construct a $\mathbf{0}''$-computable function $\gamma:\omega\times C\to\{0,1\}$, and we define
\[
 \underline{P_i}^{[\varepsilon]} := \{ x\in C : \gamma(i,x) = \varepsilon\}, \ \text{for } i\in\omega \text{ and } \varepsilon \in \{ 0,1\}.
\]

Suppose that $C = \{ c_0 < c_1 < c_2 <\dots \}$. We build the map $\gamma$ by finite extensions. At a stage $s+1$, we proceed as follows. 

\textbf{(1)} First, we consider each nonempty set $F \subseteq \{0,1,\dots,s\}$ and each map $\eta : F\to\{0,1\}$. We choose the least index $t$ such that for all $i\in F$, the value $\gamma(i, c_t)$ is not yet defined and we put $\gamma(i,c_t) := \eta(i)$ for all $i\in F$.

This ensures that $\operatorname{dom}(\gamma) = \omega \times C$ and that for each finite map $\eta \colon F \to \{ 0,1\}$, where $F$ is an \emph{arbitrary} nonempty finite subset of $\omega$, there exist infinitely many $c_t$ such that $\gamma(i,c_t) = \eta(i)$ for all $i\in F$. In verification, we refer to this fact as the \emph{pattern propagation}.

\smallskip

\textbf{(2)} Next we satisfy Eq.~(\ref{eq:free}). Notice that if $G \subseteq F$, then we always have
\[
	\mathbf{p}_F = \bigwedge_{i\in F} \mathbf{p}_i \leq \bigwedge_{i\in G} \mathbf{p}_i = \mathbf{p}_G,
\]
that is, the implication from right to left is always satisfied. Hence, we need to force the implication from left to right in Eq.~(\ref{eq:free})

For each nonempty $F\subseteq \{0,1,\dots,s\}$ and each nonempty $\mathcal{G} \subseteq \PP(\{0,1,\dots,s\}) \setminus \{\emptyset\}$ such that no $G\in\mathcal G$ is a subset of $F$, we proceed as follows. Suppose that $F = \{ b_1 < b_2 < \dots < b_k\}$.

Recall that by~\cite{DvR90} the infimum $\mathbf{p}_{F}$ contains the following `canonical' representative $\underline{Q_{F}}$:
\[
	(\underline{Q_{F}})^{[\varepsilon]} = \bigg( \prod_{i\in F} \underline{P_i} \bigg)^{[\varepsilon]} = \bigg\{ \langle y_1,y_2,\dots, y_k\rangle : y_i \in \underline{P_{b_i}}^{[\varepsilon]} \bigg\}.
\]
Similarly, we can define a canonical representative $\underline{R_{\mathcal{G}}} \in \bigvee_{G\in\mathcal{G}} \mathbf{p}_{G}$: essentially, $\underline{R_{\mathcal{G}}} = \bigoplus_{G\in \mathcal{G}} \underline{Q_{G}}$. 

Put $e := s$. If the computable function $\varphi_e$ is total (note that this could be checked with the oracle $\mathbf{0}''$), then we view $\varphi_e$ as a potential $\sm$-reduction from $\underline{Q_F}$ to $\underline{R_{\mathcal{G}}}$, and we want to diagonalize against $\varphi_e$. Thus, we proceed as follows.

We choose the least index $t$ such that the values $\gamma(i,c_t)$ are undefined for all $i \in \{0,1,\dots,s\}$. Consider the number
\[
	\delta(c_t) = \underbrace{\langle c_t,c_t,\dots,c_t\rangle}_{k\text{ times}}.
\]
We find the set $G\in \mathcal{G}$ such that the value $\varphi_e(\delta(c_t))$ belongs to the $\underline{Q_G}$-summand in the uniform join $\underline{R_{\mathcal{G}}}$. Without loss of generality, we may assume that
\[
	\varphi_e(\delta(c_t)) = \langle z_{k_0},z_{k_1},\dots,z_{k_{m}}\rangle, \text{ where } G = \{ k_0 < k_1 < \dots < k_m\}.
\]
Since $G\not\subseteq F$, we also may assume that, say, the index $k_0$ belongs to $G\setminus F$. If $z_{k_0}\not \in C$, then we can choose the values $\gamma(i,c_t)$, for $i\in F$, as some fixed $\varepsilon$, since in this case $\varphi_e(\delta(c_t))$ belongs to neither component $(\underline{R_{\mathcal{G}}})^{[\varepsilon]}$, for $\varepsilon \in \{ 0,1\}$.

Suppose that $z_{k_0}\in C$. If the value $\gamma(k_0,z_{k_0})$ is still undefined, then we put $\gamma(k_0,z_{k_0}) := 0$. We set $\gamma(\ell,c_t) := 1 - \gamma(k_0,z_{k_0})$ for all $\ell\in F$. (Note that this assignment is valid even if $c_t = z_{k_0}$, since $k_0 \not\in F$). This ensures that 
\[
	\varphi_e \colon \underline{Q_{F}} \not\smle \underline{R_{\mathcal{G}}}.
\]
Indeed, consider $\varepsilon_0 := \gamma(k_0,z_{k_0})$. Then the element $\delta(c_t)$ belongs to $\underline{Q_{F}}^{[1-\varepsilon_0]}$, while $\varphi_e(\delta(c_t))$ cannot belong to $\underline{R_{\mathcal{G}}}^{[1-\varepsilon_0]}$.

\smallskip

\textbf{(3)} We satisfy Eq.~(\ref{eq:omit}). If the function $\varphi_e$ is total, then for each $i\leq s$ we proceed as follows.
We choose the least index $t$ such that the value $\gamma(i,c_t)$ is undefined.

We find the index $j_0\neq i$ such that the value $y_0 = \varphi_e(c_t)$ belongs to the $\underline{P_{j_0}}$-th summand of the uniform join 
$\bigoplus_{j\ne i}\underline{P_j}$. Similarly to Item~(2), if $y_0 \not\in C$, then we can choose the value $\gamma(i,c_t)$ arbitrarily.

Suppose that $y_0\in C$ and that $\gamma(j_0,y_0) = \varepsilon_0$ (if $\gamma(j_0,y_0)$ was undefined before, then we just put $\gamma(j_0,y_0) = 0$). We put $\gamma(i,c_t) = 1-\varepsilon_0$. Similarly to Item~(2) one can show that now we have
\[
	\varphi_e \colon \underline{P_i} \not\smle \bigoplus_{j\neq i} \underline{P_j}.
\]
This concludes the construction.

\smallskip

\emph{Verification.} It is immediate that the constructed degrees $\mathbf{p}_i = \dsm(\underline{P_i})$ satisfy Items~(b) and~(c) of the proposition.

Item~(1) of the construction ensures that for every $i\in\omega$ and $\varepsilon \in \{ 0,1\}$, the set $\underline{P_i}^{[\varepsilon]}$ is infinite. Thus, Lemma~\ref{lem:cohesiveprime} implies that every $\mathbf{p}_i$ satisfies $\Prm(\mathbf{p}_i)$.

Towards a contradiction, assume that $\mathbf{p}_F = \mathbf{0}$ for some finite nonempty set $F$. Without loss of generality, we may assume that $0\in F$. We obtain that there is a computable set $H$ separating the sets $(\underline{Q_{F}})^{[0]} \subseteq H$ and $(\underline{Q_{F}})^{[1]} \subseteq \overline{H}$. We define a new computable set $U$ as follows: $x\in U$ if and only if $\langle x,x,\dots,x\rangle \in H$. Since every pattern propagates (including having all zeros on $F$, recall Item~(1)), there are infinitely many $x$ belonging to $C\cap U$. On the other hand, (by looking at all ones on $F$) there are infinitely many $y$ in $C \cap \overline{U}$. This contradicts the $r$-cohesiveness of the set $C$, hence we obtain Item~(a) of the proposition.

Now note that every term in a distributive lattice with $\mathbf{0}$ has a finite disjunctive normal form as a join of nonempty finite meets (or the empty join). Formula~(\ref{eq:free}) gives precisely the order comparisons between these normal forms in the free distributive lattice, and Item~(a) handles comparisons with $\mathbf{0}$. Hence, the set $\{\mathbf{0}\} \cup \{\mathbf{p}_i : i\in\omega\}$ generates (inside $\Lsm$) a free distributive lattice with zero.

Finally, if we take $\mathbf{a} = \deg_{sm}(\bigoplus_{i\in\omega} \underline{P_i})$, then $\mathbf{p}_i < \mathbf{a}$ satisfies $\Prm(\mathbf{p}_i)$. In addition,
\[
	\underline{P_i} \oplus \bigoplus_{j\neq i} \underline{P_j} \in \mathbf{a}.
\]
Hence, $\mathbf{p}_i$ belongs to $\CP_{\mathbf{a}}[\Lsm]$. 
Proposition~\ref{prop:free} is proved.
\end{proof}

\medskip

\section{Coding graphs into separation degrees}\label{sec:relations}

In this section, we show how to encode a given countable graph $G$ into $\Lsm$ using parameters.
Here an \emph{undirected graph} $G = (V,E)$ consists of a countable set of vertices $V$ and the edge relation $E$. We always assume that the relation $E$ is symmetric and irreflexive.

\begin{definition}\label{def:Gabg}
	Given $\sm$-degrees $\mathbf{a}$, $\mathbf{b}$, and $\mathbf{g}$ such that $\mathbf{b} \leq \mathbf{a}$ and $\mathbf{g} \leq \mathbf{a}$, we define an undirected graph $G_{\mathbf{a},\mathbf{b},\mathbf{g}}$ as follows.
	\begin{itemize}
		\item Its set of vertices $V$ is defined by the formula
		\[
			V_{\mathbf{a},\mathbf{b}}(x) := \CP_{\mathbf{a}}(x) \,\&\, (x\leq \mathbf{b}).
		\]
		Note that by Lemma~\ref{lem:Ja}, the set $V_{\mathbf{a},\mathbf{b}}$ is at most countable.
		
		\item The edge relation $E$ is defined by the formula
		\begin{equation*} 
			 E_{\mathbf{a},\mathbf{b},\mathbf{g}}(x,y) := V_{\mathbf{a},\mathbf{b}}(x) \,\&\, V_{\mathbf{a},\mathbf{b}}(y) \,\&\, (x\ne y) \,\&\, (x\wedge y\leq \mathbf{g}).
		\end{equation*}
	\end{itemize}
\end{definition}

\smallskip

We establish the following result about the graphs $G_{\mathbf{a},\mathbf{b},\mathbf{g}}$:

\begin{proposition}\label{prop:primegraphs}
	There exist $\sm$-degrees $\mathbf{b} \le \mathbf{a}$ such that the set $V_{\mathbf{a},\mathbf{b}}$ is countably infinite, and for every countably infinite undirected graph $G$ there is a degree $\mathbf{g} \le \mathbf{a}$ such that the graphs $G$ and $G_{\mathbf{a},\mathbf{b},\mathbf{g}}$ are isomorphic.
\end{proposition}
\begin{proof}
We take the degrees $\mathbf{p}_i = \dsm(\underline{P_i})$ constructed in
Proposition~\ref{prop:free}, and we define
\[
	\mathbf{a} = \dsm\bigg(\bigoplus_{i\in\omega} \underline{P_i} \bigg).
\] 
By Proposition~\ref{prop:free}, we have
\[
	\{\mathbf{p}_i : i\in\omega\} \subseteq \CP_{\mathbf{a}}[\Lsm].
\]
We apply Proposition~\ref{prop:fullsets}(a) to the sets $D = \CP_{\mathbf{a}}[\Lsm]$ and $S = \{\mathbf{p}_i : i\in\omega\}$.
This gives us a degree $\mathbf{b}\leq \mathbf{a}$ such that
\[
	\{ \mathbf{p}_i : i\in\omega\} = \{ \mathbf{x} \in \CP_{\mathbf{a}}[\Lsm] : \mathbf{x}\leq \mathbf{b}\}.
\]
Hence, we have $V_{\mathbf{a},\mathbf{b}}[\Lsm] = \{ \mathbf{p}_i : i\in\omega\}$.

Suppose that $G = (V,E)$ is a countably infinite undirected graph. Without loss of generality, we may assume that $V = \N$. 

For indices $i<j$, we define $\mathbf{t}_{i,j} = \mathbf{p}_i\wedge \mathbf{p}_j$. By Proposition~\ref{prop:free}(a), the elements $\mathbf{t}_{i,j}$ are nonzero. In addition, Proposition~\ref{prop:free}(b) shows that the family 
\[
	T = \{\mathbf{t}_{i,j} : i<j \} 
\]
is join-independent (recall Definition~\ref{def:join-ind}). Indeed, if $\mathbf{t}_{i,j} \leq \bigvee_{(k,\ell)\in F} \mathbf{t}_{k,\ell}$ for a nonempty finite set $F$, then Eq.~(\ref{eq:free}) implies that $\{k,\ell\} \subseteq \{i,j\}$ for some $(k,\ell) \in F$ and hence, $(k,\ell) = (i,j)$. 
Notice that, in particular, the degrees $\mathbf{t}_{i,j}$ are pairwise distinct.

Now we apply Lemma~\ref{lem:trace} to the degree $\mathbf{a}$, and the sets $T$ and 
\[
	S_E = \{ \mathbf{t}_{i,j} : i < j \text{ and } E(i,j)\}.
\]
This gives us a degree $\mathbf{g}\leq \mathbf{a}$ satisfying the condition: for $i<j$,
\[
	\mathbf{t}_{i,j} = \mathbf{p}_i \wedge \mathbf{p}_j \leq \mathbf{g}\ \Longleftrightarrow\ E(i,j).
\]
This implies that the graph $G_{\mathbf{a},\mathbf{b},\mathbf{g}}$ is isomorphic to the graph $G$.
\end{proof}

\medskip


\section{Finishing the proof of the main result}\label{sec:finish}

First, recall the following (essentially folklore) result:
\begin{proposition}[see~\cite{Hodges} and Theorem~5.4 in~\cite{ASS}]\label{prop:hodges}
	There exists a computable undirected graph $G = (V,E)$ such that $\mathcal{N} = (\N;+,\times)$ is first-order definable in $G$ without parameters. In addition, the domain of the definable copy of $\mathcal{N}$ is a subset of $V$.
\end{proposition}

\subsection{Interpreting $\TAtwo$ in $\Lsm$ with parameters}\label{subsec:param-arithm}

By applying Proposition~\ref{prop:primegraphs} to the graph $G$ from Proposition~\ref{prop:hodges} we obtain that there exist parameters $\mathbf{a}_0,\mathbf{b}_0,\mathbf{g}_0$ such that the graph $G_{\mathbf{a}_0,\mathbf{b}_0,\mathbf{g}_0}$ (from Definition~\ref{def:Gabg}) is isomorphic to $G$.

For the graph $G$ we fix graph formulas $U_*(x)$, $A_*(x,y,z)$, and $M_*(x,y,z)$ defining the domain $\N$ and the graphs of addition and multiplication inside $G$.

For arbitrary lattice parameters $\mathbf{a},\mathbf{b},\mathbf{g}$, in a natural way we translate graph formulas to lattice formulas (by using the definition of the graph $G_{\mathbf{a},\mathbf{b},\mathbf{g}}$). 
Namely, for a given graph formula $\psi(\bar x)$, the formula $\psi^G(\bar x)$ is obtained as follows: in the formula $\psi^G$ each quantifier from $\psi$ is restricted to the set $V_{\mathbf{a},\mathbf{b}}$, and each subformula of the form $E(x,y)$ is replaced by
the formula $E_{\mathbf{a},\mathbf{b},\mathbf{g}}(x,y)$. 

We define
\begin{align*}
 U(x) & := V_{\mathbf{a},\mathbf{b}}(x) \,\&\, U_*^G(x),\\
 A(x,y,z) & :=
 U(x) \,\&\, U(y) \,\&\, U(z) \,\&\, A_*^G(x,y,z),\\
 M(x,y,z) & := 
 U(x) \,\&\, U(y) \,\&\, U(z) \,\&\, M_*^G(x,y,z).
\end{align*}

Now we give formulas providing the first-order definition of $(\N,\PP(\mathbb{N}); +,\times,\in)$ inside $\Lsm$ (with parameters $\mathbf{a},\mathbf{b},\mathbf{g}$):
\begin{itemize}
	\item The set $\mathbb{N}$ is interpreted as the predicate $U(x)$.
	
	\item The formulas $A(x,y,z)$ and $M(x,y,z)$ interpret the graphs of addition and multiplication. 
	
	\item The set $\PP(\N)$ is interpreted as the set of its codes $\{ \mathbf{v} : \mathbf{v}\leq \mathbf{a}\}$.
	
	\item The relation $x\in V$ is interpreted as $\mathbf{x} \leq \mathbf{v}$ (if $x\in \mathbb{N}$ is viewed as $\mathbf{x}$ satisfying $U(\mathbf{x})$ and $V\subseteq \N$ is viewed as its code $\mathbf{v}\leq \mathbf{a}$).
	
	\item We take the quotient of the second-sort part (i.e., the interpretation of $\PP(\N)$) via the definable equivalence relation on codes of sets
	\begin{equation}\label{eq:seteq}
		 \EqSet_U(v,w) :=
 		\forall x\,(U(x)\longrightarrow(x\leq v \leftrightarrow x\leq w)).
	\end{equation}
\end{itemize}

Now put $(\mathbf{a},\mathbf{b},\mathbf{g}) := (\mathbf{a}_0,\mathbf{b}_0,\mathbf{g}_0)$. Since $G_{\mathbf{a}_0,\mathbf{b}_0,\mathbf{g}_0}$ is isomorphic to the graph $G$, our interpretation provides an isomorphic copy of $(\mathbb{N};+,\times)$. By Proposition~\ref{prop:fullsets}(a) applied to $D=U$ and arbitrary $S\subseteq U$, there exists a degree $\mathbf{c}_{S} \leq \mathbf{a}$ satisfying 
\[
	S = \{ \mathbf{x} : U(\mathbf{x}) \,\&\, (\mathbf{x} \leq \mathbf{c}_S)\}. 
\]
This implies that every set $V\in \PP(\N)$ gets its code. Similarly to the proof of Proposition~\ref{prop:fullsets}(b), now one can show that our interpretation provides an isomorphic copy of $(\N,\PP(\N);+,\times,\in)$.


\subsection{Eliminating parameters in the interpretation}

We show how to get rid of the concrete choice of parameters $\mathbf{a}_0,\mathbf{b}_0,\mathbf{g}_0$ in the interpretation given in the previous subsection.

First, we introduce ancillary lattice formulas (essentially) defining the constants $0,1 \in \mathbb{N}$ and the graph of the successor function:
\begin{align*}
 Z(x) & := U(x) \,\&\, A(x,x,x),\\
 O(x) & := U(x) \,\&\, \neg Z(x) \,\&\, M(x,x,x),\\
 S(x,y)& := U(x) \,\&\, U(y) \,\&\, \exists v\,(O(v)\& A(x,v,y)).
\end{align*}

Intuitively, the next definition gives a `lattice-coded' version for the axioms of the arithmetic $\mathrm{Z}_2$.

\begin{definition}\label{def:code}
	The first-order formula $\Code(\mathbf{a},\mathbf{b},\mathbf{g})$ asserts the following conditions. 
	\begin{enumerate}[label=\textup{(C\arabic*)},leftmargin=*]
		\item $\mathbf{b}\leq \mathbf{a}$, $\mathbf{g}\leq \mathbf{a}$, and there exists $\mathbf{z}$ such that $Z(\mathbf{z}) \,\&\, \forall w\,(Z(w)\longrightarrow w=\mathbf{z})$. In what follows, we assume that we fix this (unique) element $\mathbf{z}$.
	
		\item The relation $S$ is the graph of a total injective function acting from $U$ to $U$, and $\forall x\,\neg S(x,\mathbf{z})$.
	
		\item The relations $A$ and $M$ are graphs of total functions from $U^2$ to $U$.
	
		\item The following recursion conditions on addition and multiplication hold:
			\begin{align*}
				 &\forall x\bigl(U(x)\longrightarrow (A(x,\mathbf{z},x)\& M(x,\mathbf{z},\mathbf{z}))\bigr),\\
 				&\forall x \forall y \forall v \forall t \forall w\, \bigl(S(y,v)\& A(x,y,t)\& S(t,w) \longrightarrow A(x,v,w)\bigr),\\
				 &\forall x \forall y \forall v \forall t \forall w\,\bigl(S(y,v)\& M(x,y,t)\& A(t,x,w)
               \longrightarrow M(x,v,w)\bigr).
			\end{align*}
		\item Full induction holds:
		\begin{multline*}
 			\forall c\Bigl[c\leq \mathbf{a}\longrightarrow \Bigl(\bigl[\mathbf{z}\leq c\ \&\ \forall x\forall y\,(S(x,y)\& x\leq c\longrightarrow y\leq c)\bigr]\\
			\longrightarrow\forall x\,(U(x)\longrightarrow x\leq c)\Bigr)\Bigr].
		\end{multline*}
	\end{enumerate}
\end{definition}
Note that the tuple $(\mathbf a_0,\mathbf b_0,\mathbf g_0)$ satisfies $\Code(\mathbf{a}_0,\mathbf{b}_0,\mathbf{g}_0)$: indeed, Conditions (C1)--(C4) hold in its copy of standard arithmetic, and (C5) follows by the ordinary induction on $\mathbb N$.

We show that the formula $\Code$ allows us to eliminate the parameters.

\begin{proposition}\label{prop:standard}
Every tuple $(\mathbf{a},\mathbf{b},\mathbf{g})$ satisfying $\Code(\mathbf{a},\mathbf{b},\mathbf{g})$ in $\Lsm$ defines standard arithmetic on the set $U$. The degrees $c\leq \mathbf{a}$, modulo the equivalence relation $\EqSet_U$ from Eq.~(\ref{eq:seteq}), represent the powerset $\PP(\mathbb{N})$. Consequently, $\Code(\mathbf{a},\mathbf{b},\mathbf{g})$ implies that the interpretation of Subsection~\ref{subsec:param-arithm} applied to $(\mathbf{a},\mathbf{b},\mathbf{g})$ interprets $(\N,\PP(\N);+,\times,\in)$ in $\Lsm$.
\end{proposition}
\begin{proof}
Let $\mathbf{z}$ be the unique zero specified by Condition~(C1). We define the function $h\colon \mathbb{N} \to \Lsm$ as follows: $h(0):=\mathbf{z}$ and $h(n+1) := S(h(n))$. Condition~(C2) makes the function $h$ injective. Its range $H\subseteq U$ contains $\mathbf{z}$ and is closed under $S$. 

Proposition~\ref{prop:fullsets}(a) applied to $D=U$ and $S=H$ gives a degree $\mathbf{c}\leq \mathbf{a}$ such that 
\[
	H = \{\mathbf{x} : U(\mathbf{x}) \,\&\, \mathbf{x} \leq \mathbf{c}\}.
\] 
Condition~(C5) now implies that $U=H$.

External induction on natural numbers $m\in \mathbb{N}$, using Conditions~(C3)--(C4), gives:
\[
 A(h(n),h(m),h(n+m)) \text{ and }M(h(n),h(m),h(nm)), \text{\ \ for}\ n,m\in\omega.
\]
Therefore, the map $h$ is an isomorphism from $(\mathbb{N};+,\times)$ onto the structure $(U;A,M)$.

Finally, the map 
\[
	\theta: \mathbf{x} \mapsto \{ n : h(n)\leq \mathbf{x}\}
\] 
is a surjection from $\{ \mathbf{x} : \mathbf{x} \leq \mathbf{a}\}$ onto $\PP(\N)$. Indeed, suppose that $V\subseteq \N$. By applying Proposition~\ref{prop:fullsets}(a) to $D=U$ and the set $S = \{h(n) : n\in V\}$, we get a degree $\mathbf{x}_V \leq \mathbf{a}$ satisfying 
\[
	 \{h(n) : n\in V\} = S = \{\mathbf{y} : U(\mathbf{y}) \,\&\, (\mathbf{y} \leq \mathbf{x}_{V})\} = \{ h(n) : h(n) \leq \mathbf{x}_V\}.
\]
Hence, $\theta(\mathbf{x}_V) = V$. We conclude that the interpretation of Subsection~\ref{subsec:param-arithm} applied to the tuple $(\mathbf{a},\mathbf{b},\mathbf{g})$ interprets $(\N,\PP(\N);+,\times,\in)$.
\end{proof}

\subsection{Final steps of the proof}

First, a standard computability-theoretic argument shows that $\Th(\Lsm) \leq_m \TAtwo = \Th(\mathbb{N},\PP(\N);+,\times,\in)$ (see, e.g., Theorem~V.7.3 in~\cite{Odi-92} as an example).

Now given a sentence $\psi$ of second-order arithmetic, we effectively transform it into a lattice-theoretic sentence $\psi^{\#}(\mathbf{a},\mathbf{b},\mathbf{g})$ with parameters $\mathbf{a},\mathbf{b},\mathbf{g}$. As usual, we may assume that every atomic subformula of $\psi$ contains at most one occurrence of the signature symbols. We sketch the transformation:
\begin{itemize}
	\item addition and multiplication are replaced by the relations $A$ and $M$ respectively, 

	\item number quantifiers are translated to quantifiers over the set $U$, 
	
	\item set variables $V$ are replaced by degrees $\mathbf{v}\leq \mathbf{a}$, 
	
	\item membership $x\in V$  is replaced by $\mathbf{x}\leq \mathbf{v}$, 
	
	\item equality of sets is replaced by $\EqSet_U$. 
\end{itemize}

Then Proposition~\ref{prop:standard} implies that
\[
	 \psi \in \TAtwo\ \Leftrightarrow\  \Lsm \models\exists a \exists b \exists g \bigl( \Code(a,b,g)\,\&\, \psi^{\#}(a,b,g)\bigr).
\]
Since the translation $\psi \mapsto \psi^{\#}$ is effective, we obtain a reduction $\TAtwo \leq_m \Th(\Lsm)$. Since every first-order theory is a cylinder, we deduce that the theories $\TAtwo$ and $\Th(\Lsm)$ are recursively isomorphic.
Theorem~\ref{thm:main} is proved.


\bigskip

\section*{Declarations}

\subsection*{AI use}
The author conceived the research program and outlined the proof strategies. ChatGPT-6.2 assisted in working out the technical details of the proofs, checking intermediate arguments, and preparing the LaTeX manuscript. The author takes full responsibility for the mathematical content of the paper.

\bibliographystyle{alpha}
\bibliography{references}

\end{document}